\documentclass[11pt]{amsart}
\usepackage{epsfig}

\topskip  \DeclareMathOperator{\spanOp}{span}
\newtheorem{theorem}{Theorem}
\newtheorem{prop}[theorem]{Proposition}
\newtheorem{lemma}[theorem]{Lemma}

\newtheorem{definition}[theorem]{Definition}
\newtheorem{example}[theorem]{Example}

\def\ms{\medskip}
\def\fv{\Phi_{\mathbf{v}}}
\def\np{{\rm NP}}
\def\nf{{\mathcal N}}
\def\f{\mathbf{f}}
\def\mn{\hspace{.88pt}{\bf -}\hspace{.88pt}}

\def\A{\texttt{A}}
\def\C{\texttt{C}}
\def\G{\texttt{G}}
\def\T{\texttt{T}}
\def\Prob{\mathrm{Prob}}
\def\rr{\mathbb{R}}

\providecommand{\abs}[1]{\left\lvert#1\right\rvert}
\providecommand{\norm}[1]{\left\lVert#1\right\rVert}

\providecommand{\Z}{\mathbb{Z}}
\providecommand{\boxend}{\hspace{\stretch{1}}$\Box$\\ \ \\}

\author{Sergi Elizalde}
\address{Department of Mathematics, Dartmouth College, Hanover, NH 03755}
\email{sergi.elizalde@dartmouth.edu}
\author{Kevin Woods}
\address{Department of Mathematics, University of California, Berkeley, CA 94720}
\email{kwoods@math.berkeley.edu}

\title[The number of inference functions]{Bounds on the number of inference functions of a graphical
model}

\begin{document}
\maketitle

\begin{abstract}
Directed and undirected graphical models, also called Bayesian networks and Markov random fields, respectively, are
important statistical tools in a wide variety of fields, ranging
from computational biology to probabilistic artificial intelligence.
We give an upper bound on the number of inference functions of any
graphical model. This bound is polynomial on the size of the model,
for a fixed number of parameters, thus improving the exponential
upper bound given by Pachter and Sturmfels~\cite{PS3}. We also show
that our bound is tight up to a constant factor, by constructing a
family of hidden Markov models whose number of inference functions
agrees asymptotically with the upper bound. Finally, we apply this
bound to a model for sequence alignment that is used in
computational biology.
\end{abstract}


\section{Introduction}

Many statistical models seek, given a set of observed data, to find
the \emph{hidden} (unobserved) data which best explains these
observations. In this paper we consider graphical models (both directed and undirected), a broad
class that includes many useful models, such as hidden Markov models
(HMMs), pairwise-hidden Markov models, hidden tree models, Markov
random fields, and some language models
(background on graphical models will be given in
Section~\ref{SubSect:GraphicalModels}).
These graphical models relate the hidden and
observed data probabilistically, and a natural problem is to
determine, given a particular observation, what is the most likely
hidden data (which is called the {\em explanation}). These models
rely on parameters that are the probabilities relating the hidden
and observed data.
Any fixed values of the parameters determine a way to assign an
explanation to each possible observation. This gives us a map,
called an {\em inference function}, from observations to
explanations.

An example of an inference function is the popular ``{\it Did you
mean}" feature from \texttt{google}, which could be implemented as a
hidden Markov model, where the observed data is what we type into
the computer, and the hidden data is what we were meaning to type.
Graphical models are frequently used in these sorts of probabilistic
approaches to machine learning, pattern recognition, and artificial
intelligence (see \cite{Jen} for an introduction).

Inference functions for graphical models are also important in
computational biology \cite[Section~1.5]{ASCB}, from where we originally drew inspiration for
this paper.  For example,
consider the \emph{gene-finding functions}, which were discussed in
\cite[Section 5]{PS2}. These inference functions (corresponding to a
particular HMM) are used to identify gene structures in DNA
sequences. An observation in such a model is a sequence of
nucleotides in the alphabet $\Sigma'=\{\A,\C,\G,\T\}$, and an
explanation is a sequence of $1$'s and $0$'s which indicate whether
the particular nucleotide is in a gene or is not. We seek to use the
information in the observed data (which we can find via DNA
sequencing) to decide on the hidden information of which nucleotides
are part of genes (which is hard to figure out directly). Another
class of examples is that of sequence alignment models \cite[Section
2.2]{ASCB}. In such models, an inference function is a map from a
pair of DNA sequences to an optimal alignment of those sequences. If
we change the parameters of the model, which alignments are optimal
may change, and so the inference functions may change.

A surprising conclusion of this paper is that there cannot be
\emph{too many} different inference functions, though the parameters
may vary continuously over all possible choices. For example, in the
homogeneous binary HMM of length 5 (see
Section~\ref{SubSect:GraphicalModels} for some definitions; they are
not important at the moment), the observed data is a binary sequence
of length 5, and the explanation will also be a binary sequence of
length 5.  At first glance, there are
\[32^{32}=1\ 461\ 501\ 637\ 330\ 902\ 918\ 203\ 684\ 832\ 716\ 283\ 019\ 655\ 932\ 542\ 976\]
possible maps from observed sequences to explanations.  In fact,
Christophe Weibel has computed that only $5266$ of these possible
maps are actually inference functions \cite{Weib}.  Indeed, for an
arbitrary graphical model, the number of possible maps from observed
sequences to explanations is, at first glance, doubly exponential in
the size of the model.  The following theorem, which we call the
\emph{Few Inference Functions Theorem}, states that, if we fix the
number of parameters, the number of inference functions is actually
bounded by a polynomial in the size of the model.

\begin{theorem}[The Few Inference Functions Theorem]\label{th:fif}
Let $d$ be a fixed positive integer. Consider a graphical model with
$d$ parameters (see Definitions \ref{def:dgm} and \ref{def:ugm} for directed and undirected graphs, respectively). Let $M$ be the \emph{complexity} of the graphical model, where complexity is given by Definitions \ref{M-dgm} and \ref{M-ugm}, respectively.
Then, the number of inference functions of the model is $O(M^{d(d-1)})$.
\end{theorem}

As we shall see, the complexity of a graphical model is often linear in the number of vertices or edges of the underlying graph.

Different inference functions represent different criteria to decide
what is the most likely explanation for each observation. A bound on
the number of inference functions is important because it indicates
how badly a model may respond to changes in the parameter values
(which are generally known with very little certainty and only
guessed at). Also, the polynomial bound given in
Section~\ref{sec:the_fif_theorem} suggests that it might be feasible
to precompute all the inference functions of a given graphical
model, which would yield an efficient way to provide an explanation
for each given observation.

This paper is structured as follows. In Section~\ref{sec:prelim} we
introduce some preliminaries about graphical models and inference
functions, as well as some facts about polytopes. In
Section~\ref{sec:the_fif_theorem} we prove Theorem~\ref{th:fif}. In Section~\ref{sec:lower_bound} we prove that our upper
bound on the number of inference functions of a graphical model is
sharp, up to a constant factor, by constructing a family of HMMs
whose number of inference functions asymptotically matches the
bound. In Section~\ref{sec:seq_align} we show that the bound is also
asymptotically tight on a model for sequence alignment which is
actually used in computational biology. In particular, this bound will be quadratic on the
length of the input DNA sequences. We conclude with a few remarks
and possible directions for further research.

\section{Preliminaries}\label{sec:prelim}

\subsection{Graphical models}
\label{SubSect:GraphicalModels}

A {\em statistical model} is a family of joint probability
distributions for a collection of discrete random variables
$\mathbf{W}=(W_1,\dots,W_m)$, where each $W_i$ takes on values in
some finite state space $\Sigma_i$. A {\em graphical model} is
represented by a graph where each vertex $v_i$ corresponds to a
random variable $W_i$. The edges of the graph represent the
dependencies between the variables. There are two major classes of
graphical models depending on whether $G$ is a directed or an
undirected graph.

We start by discussing {\em directed graphical models}, also called
{\em Bayesian networks}, which are those represented by a finite
directed acyclic graph $G$. Each vertex $v_i$ has an associated
probability map
\begin{equation}\label{eqn:pi-dgm}
p_i:\left(\prod_{j:\ v_j\text{ a parent of }v_i}\Sigma_j\right)\longrightarrow
[0,1]^{\abs{\Sigma_i}}.
\end{equation}

 Given the states of each $W_j$ such that
$v_j$ is a parent of $v_i$, the probability that $v_i$ has a given
state is independent of all other vertices that are not descendants
of $v_i$, and this map $p_i$ gives that probability. In particular,
we have the equality
\begin{align*}\Prob(\mathbf{W}=\mathbf{\rho})&=\prod_{i}\Prob\left(W_i=\rho_i,\text{
given that }W_j=\rho_j\text{ for all parents $v_j$ of $v_i$}\right)\\
&=\prod_i\left(\left[p_i\left(\rho_{j_1},\ldots,\rho_{j_k}\right)\right]_{\rho_i}\right),\end{align*}
where $v_{j_i},\ldots,v_{j_k}$ are the parents of $v_i$.  Sources in
the digraph (which have no parents) are generally given the uniform
probability distribution on their states, though more general
distributions are possible. See \cite[Section 1.5]{ASCB} for general
background on graphical models.

\begin{example}
\label{Ex:HMM} The hidden Markov model (HMM) is a model with random
variables $\mathbf{X}=(X_1,\ldots,X_n)$ and
$\mathbf{Y}=(Y_1,\ldots,Y_n)$. Edges go from $X_i$ to  $X_{i+1}$ and
from $X_i$ to $Y_i$.
\end{example}

\begin{figure}[hbt]
\centering\epsfig{file=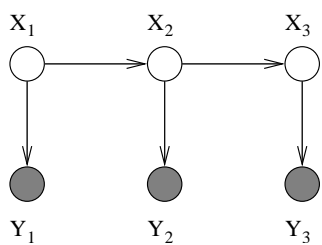,height=1.2in}
\caption{\label{fig:hmm} The graph of an HMM for $n=3$.}
\end{figure}

Generally, each $X_i$ has the same state space $\Sigma$ and each
$Y_i$ has the same state space $\Sigma'$.  An HMM is called
\emph{homogeneous} if the $p_{X_i}$, for $1\le i\le n$, are
identical and the $p_{Y_i}$ are identical.  In this case, the
$p_{X_i}$ each correspond to the same
$\abs{\Sigma}\times\abs{\Sigma}$ matrix $T=\left(t_{ij}\right)$ (the
\emph{transition matrix}) and the $p_{Y_i}$ each correspond to the
same $\abs{\Sigma}\times\abs{\Sigma'}$ matrix
$S=\left(s_{ij}\right)$ (the \emph{emission} matrix).

In the example, we have partitioned the variables into two sets.  In
general graphical models, we also have two kinds of variables:
observed variables $\mathbf{Y}=(Y_1,Y_2,\ldots,Y_n)$ and hidden
variables $\mathbf{X}=(X_1,X_2,\ldots,X_q)$. Generally, the observed
variables are the sinks of the directed graph, and the hidden
variables are the other vertices, but this does not need to be the
case. To simplify the notation, we make the assumption, which is
often the case in practice, that all the observed variables take
their values in the same finite alphabet $\Sigma'$, and that all the
hidden variables are on the finite alphabet $\Sigma$.

Notice that for given $\Sigma$ and $\Sigma'$ the homogeneous HMMs in
this example depend only on a fixed set of parameters, $t_{ij}$ and
$s_{ij}$, even as $n$ gets large.  These are the sorts of models we
are interested in.

\begin{definition}\label{def:dgm}
A \emph{directed graphical model with $d$ parameters},
$\theta_1,\ldots,\theta_d$, is a directed graphical model such that
each probability
$\left[p_i\left(\rho_{j_1},\ldots,\rho_{j_k}\right)\right]_{\rho_i}$
in (\ref{eqn:pi-dgm}) is a monomial in $\theta_1,\ldots,\theta_d$.
\end{definition}

In what follows we denote by $E$ the number of edges of the
underlying graph of a graphical model, by $n$ the number of observed
random variables, and by $q$ the number of hidden random variables.
The observations, then, are sequences in $(\Sigma')^n$ and the
explanations are sequences in $\Sigma^q$. Let $l=|\Sigma|$ and
$l'=|\Sigma'|$.

For each observation $\mathbf{\tau}$ and hidden variables
$\mathbf{h}$, $\Prob\left(\mathbf{X}=\mathbf{h}, \
\mathbf{Y}=\mathbf{\tau}\right)$ is a monomial
$f_{\mathbf{h},\mathbf{\tau}}$ in the parameters
$\theta_1,\ldots,\theta_d$. Then for each observation
$\mathbf{\tau}\in(\Sigma')^n$, the observed probability
$\Prob(\mathbf{Y}=\mathbf{\tau})$ is the sum over all hidden data
$\mathbf{h}$ of $\Prob\left(\mathbf{X}=\mathbf{h}, \
\mathbf{Y}=\mathbf{\tau}\right)$, and so
$\Prob(\mathbf{Y}=\mathbf{\tau})$ is the polynomial
$f_{\mathbf{\tau}}=\sum_{\mathbf{h}}f_{\mathbf{h},\mathbf{\tau}}$ in
the parameters $\theta_1,\ldots,\theta_d$.

\begin{definition}\label{M-dgm}
The \emph{complexity}, $M$, of a directed graphical model is
the maximum, over all $\mathbf{\tau}$, of the degree of the
polynomial $f_{\mathbf{\tau}}$.
\end{definition}

  In many graphical models, $M$ will
be a linear function of $n$, the number of observed variables.  For
example, in the homogeneous HMM, $M=E=2n-1$.

Note that we have not assumed that the appropriate probabilities sum
to 1.  It turns out that the analysis is much easier if we do not
place that restriction on our probabilities.  At the end of the
analysis, these restrictions may be added if desired (there are many
models in use, however, which never place that restriction; these
can no longer be properly called ``probabilistic'' models, but in
fact belong to a more general class of ``scoring'' models which our
analysis also encompasses).

\ms

The other class of graphical models are those that are represented
by an undirected graph. They are called {\em undirected graphical
models} and are also known as {\em Markov random fields}. As for
directed models, the vertices of the graph $G$ correspond to the random variables, but the joint probability is now
represented as a product of local functions defined on the maximal
cliques of the graph, instead of transition probabilities $p_i$
defined on the edges.

\begin{figure}[hbt]
\centering\epsfig{file=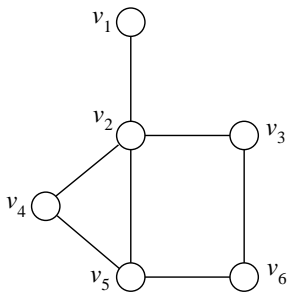,height=1.5in}
\caption{\label{fig:cliques} An undirected graph with maximal
cliques $\{v_1,v_2\}$, $\{v_2,v_3\}$, $\{v_2,v_4,v_5\}$,
$\{v_3,v_6\}$, and $\{v_5,v_6\}$. }
\end{figure}

Recall that a {\em clique} of a graph is a set of vertices with the
property that there is an edge between any two of them. A clique is
{\em maximal} if it cannot be extended to include additional
vertices without losing the property of being a clique (see
Figure~\ref{fig:cliques}).

Each maximal clique $C$ of the graph $G$ has an associated {\em
potential function}
\begin{equation}\label{eqn:pi-ugm}
\psi_C:\left(\prod_{j:\ v_j\in C} \Sigma_j\right)\longrightarrow\mathbb{R}.
\end{equation}
 Given the states $\rho_j$ of each $W_j$ such that $v_j$ is a vertex in the clique $C$,
if we denote by $\rho_C$ the vector of such states, then
$\psi_C(\rho_C)$ is a nonnegative real number. We denote by
$\mathcal{C}$ the set of all maximal cliques $C$.

Then, the joint probability distribution of all the variables $W_i$
is given by
$$\Prob(\mathbf{W}=\mathbf{\rho})=\frac{1}{Z}\prod_{C\in\mathcal{C}}\psi_C(\rho_C),$$
where $Z$ is the normalization factor
$$Z=\sum_{\mathbf{\rho}}\prod_{C\in\mathcal{C}}\psi_C(\rho_C),$$
obtained by summing over all assignments of values to the variables
$\mathbf{\rho}$.

The value of the function $\psi_C(\rho_C)$ for each possible choice
of the states $\rho_i$ is given by the parameters of the model. We
will be interested in models in which the set of parameters is
fixed, even as the size of the graph gets large.

\begin{definition}\label{def:ugm}
An \emph{undirected graphical model with $d$ parameters},
$\theta_1,\ldots,\theta_d$, is an undirected graphical model such
that each probability $\psi_C(\rho_C)$ in (\ref{eqn:pi-ugm}) is a
monomial in $\theta_1,\ldots,\theta_d$.
\end{definition}

As in the case of directed models, the variables can be partitioned
into observed variables $\mathbf{Y}=(Y_1,Y_2,\ldots,Y_n)$ (which can
be assumed to take their values in the same finite alphabet
$\Sigma'$) and hidden variables $\mathbf{X}=(X_1,X_2,\ldots,X_q)$
(which can be assumed to be on the finite alphabet $\Sigma$). For
each observation $\mathbf{\tau}$ and hidden variables $\mathbf{h}$,
$Z\cdot\Prob\left(\mathbf{X}=\mathbf{h}, \
\mathbf{Y}=\mathbf{\tau}\right)$ is a monomial
$f_{\mathbf{h},\mathbf{\tau}}$ in the parameters
$\theta_1,\ldots,\theta_d$. Then for each observation
$\mathbf{\tau}\in(\Sigma')^n$, the observed probability
$\Prob(\mathbf{Y}=\mathbf{\tau})$ is the sum over all hidden data
$\mathbf{h}$ of $\Prob\left(\mathbf{X}=\mathbf{h}, \
\mathbf{Y}=\mathbf{\tau}\right)$, and so
$Z\cdot\Prob(\mathbf{Y}=\mathbf{\tau})$ is the polynomial
$f_{\mathbf{\tau}}=\sum_{\mathbf{h}}f_{\mathbf{h},\mathbf{\tau}}$ in
the parameters $\theta_1,\ldots,\theta_d$.

\begin{definition}\label{M-ugm}
The \emph{complexity}, $M$, of an undirected graphical model is
the maximum, over all $\mathbf{\tau}$, of the degree of the
polynomial $f_{\mathbf{\tau}}$.
\end{definition}

 It is usually the case for undirected models, as in directed, that $M$ is
a linear function of $n$.

\subsection{Inference functions}

For fixed values of the parameters, the basic inference problem is
to determine, for each given observation $\mathbf{\tau}$, the value
$\mathbf{h}\in\Sigma^q$ of the hidden data that maximizes
$\Prob(\mathbf{X}=\mathbf{h}\ \big|\ \mathbf{Y}=\mathbf{\tau})$. A
solution to this optimization problem is denoted
$\mathbf{\widehat{h}}$ and is called an \emph{explanation} of the
observation $\mathbf{\tau}$. Each choice of parameter values
$(\theta_1,\theta_2,\ldots,\theta_d)$ defines an \emph{inference
function} $\mathbf{\tau}\mapsto\mathbf{\widehat{h}}$ from the set of
observations $(\Sigma')^n$ to the set of explanations $\Sigma^q$.

It is possible that there is more than one value of
$\mathbf{\widehat{h}}$ attaining the maximum of
$\Prob(\mathbf{X}=\mathbf{h}|\mathbf{Y}=\mathbf{\tau})$. In this
case, for simplicity, we will pick only one such explanation,
according to some consistent tie-breaking rule decided ahead of
time. For example, we can pick the least such $\mathbf{\widehat{h}}$
in some given total order of the set $\Sigma^q$ of hidden states.
Another alternative would be to define inference functions as maps
from $(\Sigma')^n$ to subsets of $\Sigma^q$. This would not affect
the results of this paper, so for the sake of simplicity, we
consider only inference functions as defined above.

It is interesting to observe that the total number of maps
$(\Sigma')^n\longrightarrow\Sigma^q$ is
$(l^q)^{(l')^n}=l^{q(l')^n}$, which is doubly-exponential in the
length $n$ of the observations. However, the vast majority of these maps are not
inference functions for any values of the parameters. Before our
results, the best upper bound in the literature is an exponential
bound given in \cite[Corollary 10]{PS3}. Theorem~\ref{th:fif} gives a polynomial upper bound
on the number of inference functions of a graphical model.

\subsection{Polytopes}

Here we review some facts about convex polytopes, and we introduce
some notation. Recall that a polytope is a bounded intersection of
finitely many closed halfspaces, or equivalently, the convex hull of
a finite set of points. For the basic definitions about polytopes we
refer the reader to~\cite{Zie}.

Given a polynomial $f(\theta)=\sum_{i=1}^{N} \theta_1^{a_{1,i}}
\theta_2^{a_{2,i}}\cdots \theta_d^{a_{d,i}}$, its \emph{Newton
polytope}, denoted by $\np(f)$, is defined as the convex hull in
$\mathbb{R}^d$ of the set of points
$\{(a_{1,i},a_{2,i},\ldots,a_{d,i}):i=1,\ldots,N\}$.

\begin{figure}[hbt]
\centering\epsfig{file=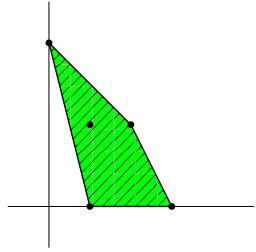,height=1.3in}
\caption{\label{fig:newton} The Newton polytope of
$f(\theta_1,\theta_2)=2\theta_1^3+3\theta_1^2\theta_2^2+\theta_1\theta_2^2+3\theta_1+5\theta_2^4$.}
\end{figure}

For example, if
$f(\theta_1,\theta_2)=2\theta_1^3+3\theta_1^2\theta_2^2+\theta_1\theta_2^2+3\theta_1+5\theta_2^4$,
then its Newton polytope $\np(f)$ is given in
Figure~\ref{fig:newton}.

Given a polytope $P \subset \rr^d$ and a vector $w \in \rr^d$, the
set of all points in $P$ at which the linear functional $\,x \mapsto
x \cdot w\,$ attains its maximum determines a {\em face} of $P$. It
is denoted
\begin{equation}
\label{faceofP}
 {\rm face}_w(P) \quad = \quad
\bigl\{\, x \in P \,\,: \,\,x \cdot w \ge y \cdot w \,\,\,\hbox{for
all} \,\, y \in P \,\bigr\}.
\end{equation}
Faces of dimension 0 (consisting of a single point) are called {\em
vertices}, and faces of dimension 1 are called {\em edges}. If $d$
is the dimension of the polytope, then faces of dimension $d-1$ are
called {\em facets}.

Let $P$ be a polytope and $F$ a face of $P$.  The {\em normal cone}
of $P$ at $F$ is
$$ N_P(F) \quad = \quad \bigl\{ w \in \rr^d \, :\, {\rm face}_w(P) =
F \,\bigr\}. $$

\begin{figure}[hbt]
\centering\epsfig{file=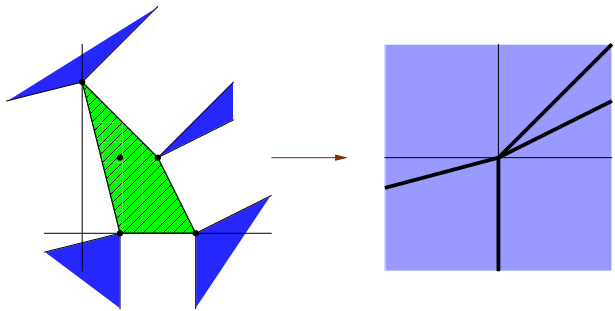,height=1.6in}
\caption{\label{fig:normal} The normal fan of a polytope.}
\end{figure}

The collection of all cones $N_P(F)$ as $F$ runs over all faces of
$P$ is denoted $\nf(P)$ and is called the \index{normal!fan} {\em
normal fan} of $P$.  Thus the normal fan $\nf(P)$ is a partition of
$\rr^d$ into cones. The cones in $\nf(P)$ are in bijection with the
faces of $P$, and if $w\in N_P(F)$ than the linear functional
$w\cdot c$ is maximized on $F$. Figure~\ref{fig:normal} shows the
normal fan of the polytope from Figure~\ref{fig:newton}.

\begin{figure}[hbt]
\centering\epsfig{file=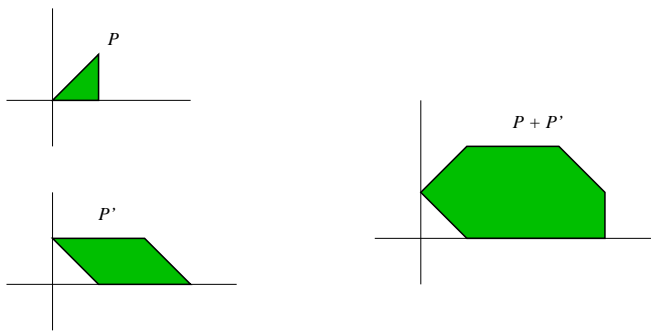,height=1.7in}
\caption{\label{fig:minksum} Two polytopes and their Minkowski sum.}
\end{figure}

The \emph{Minkowski sum} of two polytopes $P$ and $P'$ is defined as
$$P+ P':=\{\mathbf{x}+\mathbf{x'}:\mathbf{x}\in P,\,
\mathbf{x'}\in P'\}.$$ Figure~\ref{fig:minksum} shows an example in
2 dimensions. The Newton polytope of the map
$\f:\mathbb{R}^d\longrightarrow\mathbb{R}^{(l')^n}$ is defined as
the Minkowski sum of the individual Newton polytopes of its
coordinates, namely
$\np(\f):=\sum_{\mathbf{\tau}\in(\Sigma')^n}\np(f_\mathbf{\tau})$.

The \index{common refinement} \emph{common refinement} of two or
more normal fans is the collection of cones obtained as the
intersection of a cone from each of the individual fans. For
polytopes $P_1,P_2,\ldots,P_k$, the common refinement of their
normal fans is denoted $\nf(P_1)\wedge\cdots\wedge\nf(P_k)$. The
following lemma states the well-known fact that the normal fan of a
Minkowski sum of polytopes is the common refinement of their
individual fans (see \cite[Proposition 7.12]{Zie} or \cite[Lemma
2.1.5]{GS}):
\begin{lemma}\label{lem:sum_refinement}
$\nf(P_1+\cdots+ P_k)=\nf(P_1)\wedge\cdots\wedge\nf(P_k)$.
\end{lemma}

We finish with a result of Gritzmann and Sturmfels that will be
useful later.
It gives a bound on the number of vertices of a Minkowski sum of
polytopes.

\begin{theorem}[\cite{GS}]\label{thm:non-parallel}
Let $P_1,P_2,\ldots,P_k$ be polytopes in $\mathbb{R}^d$, and let $m$
denote the number of non-parallel edges of $P_1,\dots,P_k$. Then the
number of vertices of $P_1+\cdots+ P_k$ is at most
$$2 \sum_{j=0}^{d-1} \binom{m-1}{j}.$$
\end{theorem}
Note that this bound is independent of the number $k$ of polytopes.

\section{An upper bound on the number of inference functions}\label{sec:the_fif_theorem}

For fixed parameters, the inference problem of finding the
explanation $\mathbf{\widehat{h}}$ that maximizes
$\Prob(\mathbf{X}=\mathbf{h}|\mathbf{Y}=\mathbf{\tau})$ is
equivalent to identifying the monomial
$f_{\mathbf{\widehat{h}},\mathbf{\tau}}= \theta_1^{a_{1,i}}
\theta_2^{a_{2,i}}\cdots \theta_d^{a_{d,i}}$ of $f_\mathbf{\tau}$
with maximum value. Since the logarithm is a monotonically
increasing function, the desired monomial also maximizes the
quantity
\begin{eqnarray}\nonumber \log(\theta_1^{a_{1,i}}
\theta_2^{a_{2,i}}\cdots \theta_d^{a_{d,i}})&=&
a_{1,i}\log(\theta_1)+a_{2,i}\log(\theta_2)+\cdots+a_{d,i}\log(\theta_d)\\
&=&a_{1,i}v_1+a_{2,i}v_2+\cdots+a_{d,i}v_d,\nonumber\end{eqnarray}
where we replace $\log(\theta_i)$ with $v_i$. This is equivalent to
the fact that the corresponding point
$(a_{1,i},a_{2,i},\ldots,a_{d,i})$ maximizes the linear expression
$v_1x_1+\cdots+v_dx_d$ on the Newton polytope
$\np(f_\mathbf{\tau})$. Thus, the inference problem for fixed
parameters becomes a linear programming problem.


Each choice of the parameters
$\theta=(\theta_1,\theta_2,\ldots,\theta_d)$ determines an inference
function. If $\mathbf{v}=(v_1,v_2,\ldots,v_d)$ is the vector in
$\mathbb{R}^d$ with coordinates $v_i=\log(\theta_i)$, then we denote
the corresponding inference function by
$$\fv:(\Sigma')^n\longrightarrow\Sigma^q.$$
For each observation $\mathbf{\tau}\in(\Sigma')^n$, its explanation
$\fv(\mathbf{\tau})$ is given by the vertex of
$\np(f_\mathbf{\tau})$ that is maximal in the direction of the
vector $\mathbf{v}$. Note that for certain values of the parameters
(if $\mathbf{v}$ is perpendicular to a positive-dimensional face of
$\np(f_\mathbf{\tau})$) there may be more than one vertex attaining
the maximum. It is also possible that a single point
$(a_{1,i},a_{2,i},\ldots,a_{d,i})$ in the polytope corresponds to
several different values of the hidden data. In both cases, we pick
the explanation according to the tie-breaking rule determined ahead
of time. This simplification does not affect the asymptotic number
of inference functions.

Different values of $\theta$ yield different directions
$\mathbf{v}$, which can result in distinct inference functions. We
are interested in bounding the number of different inference
functions that a graphical model can have. Theorem~\ref{th:fif} gives an
upper bound which is polynomial in the size of the graphical model.
In other words, extremely few of the $l^{q(l')^n}$ functions
$(\Sigma')^n\longrightarrow\Sigma^q$ are actually inference functions.

We use the notation $f(n)\in O(g(n))$ to indicate that
$\limsup_{n\rightarrow\infty}|f(n)/g(n)|<\infty$. Similarly $f(n)\in
\Omega(g(n))$ means that
$\liminf_{n\rightarrow\infty}|f(n)/g(n)|>0$, and $f(n)\in
\Theta(g(n))$ denotes that $f(n)$ belongs to both $O(g(n))$ and
$\Omega(g(n))$.

Before proving Theorem~\ref{th:fif}, observe that usually $M$, the complexity of the graphical model, is linear in
$n$. For example, in the case of directed models,
consider the common situation where $M$ is bounded
by $E$, the number of edges of the underlying graph (this happens
when each edge ``contributes'' at most degree 1 to the monomials
$f_{\mathbf{h},\mathbf{\tau}}$, as in the homogeneous HMM). In most
graphical models of interest, $E$ is a linear function of $n$, so
the bound becomes $O(n^{d(d-1)})$. For example, the homogeneous HMM
has $M=E=2n-1$.

In the case of undirected models, if each $\psi_C(\rho_C)$ is a
parameter of the model, then $f_{\mathbf{h},\mathbf{\tau}}=Z\cdot
\Prob\left(\mathbf{X}=\mathbf{h}, \ \mathbf{Y}=\mathbf{\tau}\right)$
is a product of potential functions for each maximal clique of the
graph, so $M$ is bounded by the number of maximal cliques, which in
many cases is also a linear function of the number of vertices of
the graph. For example, this is the situation in language models
where each word depends on a fixed number of previous words in the
sentence.

\begin{proof}
In the first part of the proof we will reduce the problem of
counting inference functions to the enumeration of the vertices of a
certain polytope. We have seen that an inference function is
specified by a choice of the parameters, which is equivalent to
choosing a vector $\mathbf{v}\in\mathbb{R}^d$. The function is
denoted $\fv:(\Sigma')^n\longrightarrow\Sigma^q$, and the
explanation $\fv(\mathbf{\tau})$ of a given observation
$\mathbf{\tau}$ is determined by the vertex of
$\np(f_\mathbf{\tau})$ that is maximal in the direction of
$\mathbf{v}$. Thus, cones of the normal fan
$\nf(\np(f_\mathbf{\tau}))$ correspond to sets of vectors
$\mathbf{v}$ that give rise to the same explanation for the
observation $\mathbf{\tau}$. Non-maximal cones (i.e., those
contained in another cone of higher dimension) correspond to
directions $\mathbf{v}$ for which more than one vertex is maximal.
Since ties are broken using a consistent rule, we disregard this
case for simplicity. Thus, in what follows we consider only maximal
cones of the normal fan.

Let $\mathbf{v'}=(v'_1,v'_2,\ldots,v'_d)$ be another vector
corresponding to a different choice of parameters (see Figure
\ref{fig:inf_fct}). By the above reasoning,
$\fv(\mathbf{\tau})=\Phi_{\mathbf{v'}}(\mathbf{\tau})$ if and only
if $\mathbf{v}$ and $\mathbf{v'}$ belong to the same cone of
$\nf(\np(f_\mathbf{\tau}))$. Thus, $\fv$ and $\Phi_{\mathbf{v'}}$
are the same inference function if and only if $\mathbf{v}$ and
$\mathbf{v'}$ belong to the same cone of $\nf(\np(f_\mathbf{\tau}))$
for all observations $\mathbf{\tau}\in(\Sigma')^n$. Consider the
common refinement of all these normal fans,
$\bigwedge_{\mathbf{\tau}\in(\Sigma')^n}\nf(\np(f_\mathbf{\tau}))$.
Then, $\fv$ and $\Phi_{\mathbf{v'}}$ are the same function exactly
when $\mathbf{v}$ and $\mathbf{v'}$ lie in the same cone of this
common refinement.

\begin{figure}[h]
\begin{center} \includegraphics[width=0.6\textwidth]{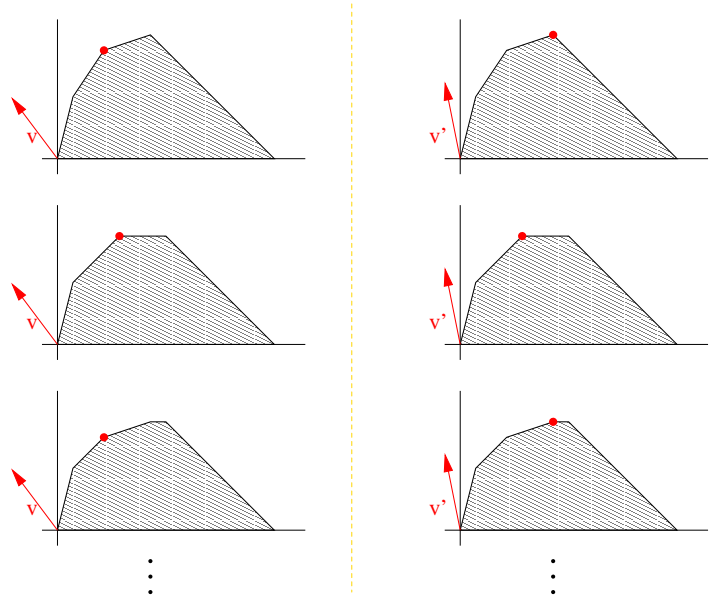}
\caption{Two different inference functions, $\fv$ (left column) and
$\Phi_{\mathbf{v'}}$ (right column). Each row corresponds to a
different observation. The respective explanations are given by the
marked vertices in each Newton polytope. \label{fig:inf_fct}}
\end{center} \end{figure}

This implies that the number of inference functions equals the
number of cones in
\[\bigwedge_{\mathbf{\tau}\in(\Sigma')^n}\nf(\np(f_\mathbf{\tau})).\]
By Lemma~\ref{lem:sum_refinement}, this common refinement is the
normal fan of
$\np(\f)=\sum_{\mathbf{\tau}\in(\Sigma')^n}\np(f_\mathbf{\tau})$,
the Minkowski sum of the polytopes $\np(f_\mathbf{\tau})$ for all
observations $\mathbf{\tau}$. It follows that enumerating inference
functions is equivalent to counting vertices of $\np(\f)$. In the
remaining part of the proof we give an upper bound on the number of
vertices of $\np(\f)$.

Note that for each $\mathbf{\tau}$, the polytope
$\np(f_\mathbf{\tau})$ is contained in the hypercube $[0,M]^d$,
since by definition of $M$, each parameter $\theta_i$ appears in
$f_\mathbf{\tau}$ with exponent at most $M$. Also, the vertices of
$\np(f_\mathbf{\tau})$ have integral coordinates, because they are
exponent vectors. Polytopes whose vertices have integral coordinates
are called \index{polytope!lattice} \emph{lattice polytopes}. It
follows that the edges of $\np(f_\mathbf{\tau})$ are given by
vectors where each coordinate is an integer between $-M$ and $M$.
There are only $(2M+1)^d$ such vectors, so this is an upper bound on
the number of different directions that the edges of the polytopes
$\np(f_\mathbf{\tau})$ can have.

This property of the Newton polytopes of the coordinates of the
model will allow us to give an upper bound on the number of vertices
of their Minkowski sum $\np(\f)$. The last ingredient that we need
is Theorem~\ref{thm:non-parallel}. In our case we have a sum of
polytopes $\np(f_\mathbf{\tau})$, one for each observation
$\mathbf{\tau}\in(\Sigma')^n$, having at most $(2M+1)^d$
non-parallel edges in total. Hence, by
Theorem~\ref{thm:non-parallel}, the number of vertices of $\np(\f)$
is at most $$2 \sum_{j=0}^{d-1} \binom{(2M+1)^d-1}{j}.$$ As $M$ goes
to infinity, the dominant term of this expression is
$$\frac{2^{d^2-d+1}}{(d-1)!}\ M^{d(d-1)}.$$
Thus, we get an $O(M^{d(d-1)})$ upper bound on the number of
inference functions of the graphical model.
\end{proof}

In the next section we will show that the bound given in
Theorem~\ref{th:fif} is tight up to a constant factor.

\section{A lower bound}\label{sec:lower_bound}

As before, we fix $d$, the number of parameters in our model. The
Few Inferences Function Theorem tells us that the number of
inference functions is bounded from above by some function
$cM^{d(d-1)}$, where $c$ is a constant (depending only on $d$) and
$M$ is the complexity of the model. Here we show that that bound is
tight up to a constant, by constructing a family of graphical models
whose number of inference functions is at least $c'M^{d(d-1)}$,
where $c'$ is another constant.  In fact, we will construct a family
of hidden Markov models with this property.  To be precise, we have
the following theorem.

\begin{theorem} \label{Thm:LowerBound} Fix $d$.  There is a constant $c'=c'(d)$ such that, given $n\in\Z_+$,
there exists an HMM of length $n$, with $d$ parameters, $4d+4$
hidden states, and $2$ observed states, such that there are at least
$c'n^{d(d-1)}$ distinct inference functions. (For this HMM, $M$ is a
linear function of $n$, so this also gives us the lower bound in
terms of $M$).
\end{theorem}

In Section \ref{SubSect:ProofLowerBound} we prove Theorem
\ref{Thm:LowerBound}. This proof requires several lemmas that we
will meet along the way, and these lemmas will be proved in Section
\ref{SubSect:LemmasLowerBound}.  Lemma \ref{Lemma:FullSet}, which is
interesting in its own right as a statement in the geometry of
numbers is proved in~\cite{EW2}.

\subsection{Proof of Theorem \ref{Thm:LowerBound}}
\label{SubSect:ProofLowerBound} Given $n$, we first construct the
appropriate HMM, $\mathcal{M}_n$, using the following lemma.

\begin{lemma}
\label{Lemma:anHMM} Given $n\in\Z_+$, there is an HMM,
$\mathcal{M}_n$, of length $n$, with $d$ parameters, $4d+4$ hidden states, and 2 observed states,
such that for any
    $a\in\Z^d_+$ with $\sum_i a_i< n$, there is an observed sequence which has one explanation
    if
    \[a_1\log(\theta_1)+\cdots+a_d\log(\theta_d)>0\]
    and another explanation if
    $a_1\log(\theta_1)+\cdots+a_d\log(\theta_d)<0$.
\end{lemma}

This means that, for the HMM $\mathcal{M}_n$, the decomposition of (log-)parameter space into
inference cones includes all of the hyperplanes $\{x:\ \langle
a,x\rangle=0\}$ such that  $a\in\Z^d_+$ with $\sum_i a_i< n$.  Call
the arrangement of these hyperplanes $\mathcal{H}_n.$ It suffices to
show that the arrangement $\mathcal{H}_n$ consists of at least
$c'n^{d(d-1)}$ chambers (full dimensional cones determined by the
arrangement). There are $c_1n^{d}$ ways to choose one of the
hyperplanes from $\mathcal{H}_n$, for some constant $c_1$. Therefore
there are $c_1^{d-1}n^{d(d-1)}$ ways to choose $d-1$ of the
hyperplanes; their intersection is, in general, a 1-dimensional face of
$\mathcal{H}_n$ (that is, the intersection is a ray which is an
extreme ray for the cones it is contained in). It is quite possible
that two different ways of choosing $d-1$ hyperplanes give the same
extreme ray.  The following lemma says that some constant fraction
of these choices of extreme rays are actually distinct.

\begin{lemma}
\label{Lemma:NumExtremeRays} Fix $d$. Given $n$, let $\mathcal{H}_n$
be the hyperplane arrangement consisting of the hyperplanes of the
form $\{x:\ \langle a,x\rangle=0\}$ with $a\in\Z^d_+$ and $\sum_i
a_i< n$. Then the number of 1-dimensional faces of $\mathcal{H}_n$
is $c_2n^{d(d-1)}$, for some constant $c_2$.
\end{lemma}

Each chamber will have a number of these extreme rays on its
boundary.  The following lemma gives a constant bound on this
number.

\begin{lemma}
\label{Lemma:ExtremeRaysPerChamber} Fix $d$. Given $n$, define
$\mathcal{H}_n$ as above. Each chamber of $\mathcal{H}_n$ has at
most $2^{d(d-1)}$ extreme rays.
\end{lemma}

Conversely, each ray is an extreme ray for at least $1$ chamber.
Therefore there are at least $\frac{c_2}{2^{d(d-1)}}n^{d(d-1)}$
    chambers, and Theorem \ref{Thm:LowerBound} is proved. \boxend

In proving Lemma \ref{Lemma:NumExtremeRays}, we will need one more
lemma. This lemma is interesting in its own right as a probabilistic
statement about integer lattices, and so is proved in a companion paper \cite{EW2}.
Given a set $S\subset\Z^d$ of
integer vectors, $\spanOp_{\rr}(S)$ is a linear subspace of $\rr^d$
and $\spanOp_{\rr}(S)\cap\Z^d$ is a sublattice of $\Z^d$.  We say
that $S$ is \emph{primitive} if $S$ is a $\Z$-basis for the
lattice $\spanOp_{\rr}(S)\cap\Z^d$.  Equivalently, a set $S$ is primitive if and only if it may be extended to a $\Z$-basis of all of $\Z^d$ (see \cite{Lek}).

We imagine picking each vector
in $S$ uniformly at random from some large box in $\rr^d$. As the
size of the box approaches infinity, the following lemma will tell
us that the probability that $S$ is primitive approaches
\[\frac{1}{\zeta(d)\zeta(d-1)\cdots\zeta(d-m+1)},\]
where $\abs{S}=m$ and $\zeta(a)$ is the Riemann Zeta function
$\sum_{i=1}^{\infty}\frac{1}{i^a}$.

\begin{lemma}[from \cite{EW2}]
\label{Lemma:FullSet} Let $d$ and $m$ be given, with $m<d$.  For
$n\in\Z_+$, $1\le k\le m$, and $1\le i\le d$, let $b_{n,k,i}\in\Z$.
For a given $n$, choose integers $s_{ki}$ uniformly (and
independently) at random from the set $b_{n,k,i}\le s_{ki}\le
b_{n,k,i}+n$.  Let $s_k=(s_{k1},\ldots,s_{kd})$ and let
$S=\{s_1,s_2,\ldots,s_m\}$.

If  $\abs{b_{n,k,i}}$ is bounded by a polynomial in $n$, then, as $n$ approaches infinity, the
probability that $S$ is a primitive set approaches
\[\frac{1}{\zeta(d)\zeta(d-1)\cdots\zeta(d-m+1)},\]
where $\zeta(a)$ is the Riemann Zeta function
$\sum_{i=1}^{\infty}\frac{1}{i^a}$.
\end{lemma}

When $m=1$, this lemma gives the probability that a $d$-tuple of
integers are relatively prime as $\frac{1}{\zeta(d)}$. For $m=1,
d=2$, this is a classic result in number theory (see \cite{Apo}),
and for $m=1, d>2$, this was proven in \cite{Nym}. Note also that,
if $m=d$ and we choose $S$ of size $m$, then the probability that
$S$ is primitive (i.e., that it is a basis for $\Z^d$) approaches zero.
This agrees with the lemma in the sense that we would expect the
probability to be
\[\frac{1}{\zeta(d)\zeta(d-1)\cdots\zeta(1)},\]
but $\zeta(1)$ does not converge.

\subsection{Proofs of Lemmas}
\label{SubSect:LemmasLowerBound}

\begin{proof}[Proof of Lemma \ref{Lemma:anHMM}]

Given $d$ and $n$, define a length $n$ HMM with parameters
$\theta_1,...,\theta_d$, as follows.  The observed states will be S
and C (for ``start of block,'' and ``continuing block,''
respectively). The hidden states will be $s_i$, $s'_i$, $c_i$, and
$c'_i$, for $1\le i\le d+1$ (think of $s_i$ and $s'_i$ as ``start of
the $i$th block'' and $c_i$ and $c'_i$ as ``continuing the $i$th
block'').

Here is the idea of what we want this HMM to do: if the observed
sequence has S's in position 1, $a_1+1$, $a_1+a_2+1$, $\ldots$, and
$a_1+\cdots+a_d+1$ and C's elsewhere, then there will be only two
possibilities for the sequence of hidden states, either
\[t=s_1\underbrace{c_1\cdots c_1}_{a_1-1}s_2\underbrace{c_2\cdots
c_2}_{a_2-1}\cdots s_d\underbrace{c_d\cdots
c_d}_{a_d-1}s_{d+1}\underbrace{c_{d+1}\cdots
c_{d+1}}_{n-a_1-\cdots-a_d-1}\] or
\[t'=s'_1\underbrace{c'_1\cdots c'_1}_{a_1-1}s'_2\underbrace{c'_2\cdots
c'_2}_{a_2-1}\cdots s'_d\underbrace{c'_d\cdots
c'_d}_{a_d-1}s'_{d+1}\underbrace{c'_{d+1}\cdots
c'_{d+1}}_{n-a_1-\cdots-a_d-1}.\] We will also make sure that $t$
has a priori probability $\theta_1^{a_1}\cdots\theta_d^{a_d}$ and
$t'$ has a priori probability 1. Then $t$ is the explanation if
$a_1\log(\theta_1)+\cdots+a_d\log(\theta_d)>0$ and $t'$ is the
explanation if $a_1\log(\theta_1)+\cdots+a_d\log(\theta_d)<0$.
 Remember that we are not constraining our probability sums to be 1.
A very similar HMM could be constructed that obeys that constraint,
if desired. To simplify notation it will be more convenient to
treat the transition probabilities as parameters that do not
necessarily sum to one at each vertex, even if this forces us to use
the term ``probability" somewhat loosely.

Here is how we set up the transitions/emmisions. Let $s_i$ and
$s'_i$, for $1\le i\le d+1$, all emit S with probability 1 and C
with probability 0. Let $c_i$ and $c'_i$ emit C with probability 1
and S with probability 0. Let $s_i$, for $1\le i\le d$, transition
to $c_i$ with probability $\theta_i$ and transition to everything
else with probability 0. Let $s_{d+1}$ transition to $c_{d+1}$ with
probability $1$ and to everything else with probability 0. Let
$s'_i$, for $1\le i\le d+1$, transition to $c'_i$ with probability
$1$ and to everything else with probability 0. Let $c_i$, for $1\le
i\le d$, transition to $c_i$ with probability $\theta_i$, to
$s_{i+1}$ with probability $\theta_{i}$, and to everything else with
probability 0. Let $c_{d+1}$ transition to $c_{d+1}$ with
probability $1$, and to everything else with probability 0. Let
$c'_i$, for $1\le i\le d$ transition to $c'_i$ with probability $1$,
to $s_{i+1}$ with probability $1$, and to everything else with
probability 0.  Let $c'_{d+1}$ transition to $c'_{d+1}$ with
probability $1$ and to everything else with probability 0.

Starting with the uniform probability distribution on the first
hidden state, this does exactly what we want it to: given the
correct observed sequence, $t$ and $t'$ are the only explanations,
with the correct probabilities.\end{proof}

\begin{proof}[Proof of Lemma \ref{Lemma:NumExtremeRays}]
We are going to pick $d-1$ vectors $a^{(1)},\ldots,a^{(d-1)}$ which
correspond to the $d-1$ hyperplanes $\{x:\ \langle
a^{(i)},x\rangle=0\}$ that will intersect to give us extreme rays
of our chambers.  We will restrict the region from which we
pick each $a^{(i)}\in\Z^d$.  Let
\[b^{(i)}=(1,1,\ldots,1)-\frac{1}{2}e_i,\]
for $1\le i\le d-1$, where $e_i$ is the $i$th standard basis vector. Let
$s=\frac{1}{4d+4}$.  For $1\le i\le d-1$, we will choose
$a^{(i)}\in\Z^d$ such that
\begin{equation}\label{eqn:C1}
\norm{\frac{n}{d}b^{(i)}-a^{(i)}}_{\infty}<\frac{n}{d}s.
\end{equation}
Note that $\sum_j a^{(i)}_j <n$, so there are observed sequences
which give us the hyperplanes $\{x:\ \langle a^{(i)},x\rangle=0\}$.
Note also that there are $(\frac{2s}{d})^{d(d-1)}n^{d(d-1)}$ choices
for the $(d-1)$-tuple of vectors $(a^{(1)},\ldots,a^{(d-1)})$.  To
prove this lemma, we must then show that a positive fraction of
these actually give rise to distinct extreme rays
$\bigcap_{i=1}^{d-1}\{x:\ \langle a^{(i)},x\rangle=0\}$.

First,  we imagine choosing the $a^{(i)}$ uniformly at random in the range given by (\ref{eqn:C1}), this probability distribution meets the condition in the
statement of Lemma \ref{Lemma:FullSet}, as $n$ approaches infinity.
 Therefore, there is a positive probability that
 \begin{equation}\label{eqn:C2}
 \{a^{(i)}:\ 1\le i \le d-1\}\text{ form
a basis for the lattice }\Z^d\cap\spanOp\{a^{(i)}:\ 1\le i\le d-1\},
\end{equation} and
this probability approaches
\[\frac{1}{\zeta(d)\zeta(d-1)\cdots\zeta(2)}.\]

Second, we look at all choices of $a^{(i)}\in\Z^d$ such that
(\ref{eqn:C1}) and (\ref{eqn:C2}) hold.  There are $c_2n^{d(d-1)}$
of these, for some constant $c_2$.  We claim that these give
distinct extreme rays $\bigcap_{i=1}^{d-1}\{x:\ \langle
a^{(i)},x\rangle=0\}$.  Indeed, say that $a^{(i)}$ and $c^{(i)}$ are both chosen such that (\ref{eqn:C1}) and (\ref{eqn:C2}) hold and such that
\[\bigcap_{i=1}^{d-1}\{x:\ \langle
a^{(i)},x\rangle=0\}=\bigcap_{i=1}^{d-1}\{x:\ \langle
c^{(i)},x\rangle=0\}.\]
We will argue that $a^{(i)}$ and $c^{(i)}$ are ``so close''
that they must actually be the same.

Let $j$, for $1\le j\le d-1$ be given.  We will prove that $a^{(j)}=c^{(j)}$. Since
\[\bigcap_{i=1}^{d-1}\{x:\ \langle
a^{(i)},x\rangle=0\}\subset\{x:\ \langle c^{(j)},x\rangle=0\},\] we
know that $c^{(j)}$ is in $\spanOp\{a^{(i)}:\ 1\le i\le d-1\}$, and therefore
\[c^{(j)}\in\Z^d\cap\spanOp\{a^{(i)}:\ 1\le i\le d\}.\]  Let $g=c^{(j)}-a^{(j)}$.  Then
\[\norm{g}_{\infty}<2\frac{n}{d}s,\] by Condition (\ref{eqn:C1}) for $a^{(i)}$ and $c^{(i)}$, and
\[g=\alpha_1a^{(1)}+\cdots+\alpha_{d-1}a^{(d-1)},\]
for some $\alpha_i\in\Z$, by Condition (\ref{eqn:C2}) for $a^{(i)}$.  We must show that $g=0$.  By reordering indices and possibly
considering $-g$, we may assume that $\alpha_1,\ldots,\alpha_k\ge
0$, for some $k$, $\alpha_{k+1},\ldots,\alpha_{d-1}\le 0$, and
$\abs{\alpha_1}$ is maximal over all $\abs{\alpha_i}$, $1\le i\le
d-1$.

Examining the first coordinate  of $g$, we have that
\begin{align*} -2\frac{n}{d}s&<g_1\\
&=\alpha_1a^{(1)}_1+\cdots+\alpha_{d-1}a^{(d-1)}_1\\
&<\alpha_1\frac{n}{d}(b^{(1)}_1+s)+\cdots+\alpha_k\frac{n}{d}(b^{(k)}_1+s)+
\alpha_{k+1}\frac{n}{d}(b^{(k+1)}_1-s)+\cdots+\alpha_{d-1}\frac{n}{d}(b^{(d-1)}_1-s)\\
&=\frac{n}{d}\big[\alpha_1+\cdots+\alpha_{d-1}-\frac{1}{2}\alpha_1+s(\abs{\alpha_1}+\cdots+\abs{\alpha_{d-1}})\big]
\ \ \ \ \text{ (using 
$b^{(i)}=(1,\ldots,1)-\frac{1}{2}e_i$)}
\\
&\le
\frac{n}{d}\big[\alpha_1+\cdots+\alpha_{d-1}-\frac{1}{2}\alpha_1+(d-1)s\alpha_1\big].\end{align*}
Negating and dividing by $\frac{n}{d}$,
\begin{equation}\label{eq:1}-(\alpha_1+\cdots+\alpha_{d-1})+\frac{1}{2}\alpha_1-(d-1)s\alpha_1<2s.\end{equation}
Similarly, examining the $(k+1)$-st coordinate of $g$, we have
\begin{align*}
2\frac{n}{d}s&>g_{k+1}\\
&=\alpha_1a^{(1)}_{k+1}+\cdots+\alpha_{d-1}a^{(d-1)}_{k+1}\\
&>\alpha_1\frac{n}{d}(b^{(1)}_{k+1}-s)+\cdots+\alpha_k\frac{n}{d}(b^{(k)}_{k+1}-s)+
\alpha_{k+1}\frac{n}{d}(b^{(k+1)}_{k+1}+s)+\cdots+\alpha_{d-1}\frac{n}{d}(b^{(d-1)}_{k+1}+s)\\
&=\frac{n}{d}\big[\alpha_1+\cdots+\alpha_{d-1}-\frac{1}{2}\alpha_{k+1}-s(\abs{\alpha_1}+\cdots+\abs{\alpha_{d-1}})\big]\\
&\ge
\frac{n}{d}\big[\alpha_1+\cdots+\alpha_{d-1}-\frac{1}{2}\alpha_{k+1}-(d-1)s\alpha_1\big],
\end{align*}
and so
\begin{equation}\label{eq:2}(\alpha_1+\cdots+\alpha_{d-1})-\frac{1}{2}\alpha_{k+1}-(d-1)s\alpha_1<2s.\end{equation}
Adding the equations (\ref{eq:1}) and (\ref{eq:2}),
\[\frac{1}{2}\alpha_1-\frac{1}{2}\alpha_{k+1}-2(d-1)s\alpha_1<4s,\]
and so, since $s=\frac{1}{4d+4}$,
\[\frac{1}{d+1}\alpha_1-\frac{1}{2}\alpha_{k+1}<\frac{1}{d+1}.\]
Therefore, since $\alpha_{k+1}\le 0$, we have that $\alpha_1<1$ and
so $\alpha_1=0$. Since $\abs{\alpha_1}$ was maximal over all
$\abs{\alpha_i}$, we have that $g=0$.  Therefore $a^{(j)}=c^{(j)}$, and the lemma follows.
\end{proof}

\begin{proof}[Proof of Lemma \ref{Lemma:ExtremeRaysPerChamber}]
Suppose $N>2^{d(d-1)}$, and suppose $a^{(i,j)},$ for $1\le i\le N$
and $1\le j\le d-1$, are such that $a^{(i,j)}\in\Z_+^d$,
$\sum_{k=1}^d a^{(i,j)}_k<n$, and the $N$ rays
\[r^{(i)}=\bigcap_{j=1}^{d-1}\{x:\ \langle a^{(i,j)},x\rangle=0\}\]
are the extreme rays for some chamber.  Then, since $N>2^{d(d-1)}$,
there are some $i$ and $i'$ such that
\[a^{(i,j)}_k\equiv a^{(i',j)}_k \mod{2},\]
for $1\le j\le d-1$ and $1\le k\le d$ (i.e., all of the coordinates
in all of the vectors have the same parity).  Then let
\[c^{(j)}=\frac{a^{(i,j)}+a^{(i',j)}}{2},\]
for $1\le j\le d-1$.  Then $c^{(j)}\in\Z_+^d$ and $\sum_{k=1}^d
c^{(j)}_k<n$, and the ray
\[r=\bigcap_{j=1}^{d-1}\{x:\ \langle c^{(j)},x\rangle=0\}=\frac{r^{(i)}+r^{(i')}}{2}\]
is in the chamber, which is a contradiction.
\end{proof}

\section{Inference functions for sequence alignment}
\label{sec:seq_align}

In this section we give an application of Theorem~\ref{th:fif} to a
basic model for sequence alignment. Sequence alignment is one of the
most frequently used techniques in determining the similarity
between biological sequences. In the standard instance of the
sequence alignment problem, we are given two sequences (usually DNA
or protein sequences) that have evolved from a common ancestor via a
series of mutations, insertions and deletions. The goal is to find
the best alignment between the two sequences. The definition of
``best" here depends on the choice of scoring scheme, and there is
often disagreement about the correct choice. In \emph{parametric
sequence alignment}, this problem is circumvented by instead
computing the optimal alignment as a function of \emph{variable}
scores. Here we consider one such scheme, in which all matches are
equally rewarded, all mismatches are equally penalized and all
spaces are equally penalized. Efficient parametric sequence
alignment algorithms are known (see for
example~\cite[Chapter~7]{ASCB}). Here we are concerned with the
different inference functions that can arise when the parameters
vary. For a detailed treatment on the subject of sequence alignment,
we refer the reader to \cite{Gus}.

Given two strings $\sigma^1$ and $\sigma^2$ of lengths $n_1$ and
$n_2$ respectively, an \emph{alignment} is a pair of equal length
strings $(\mu^1,\mu^2)$ obtained from $\sigma^1,\sigma^2$ by
inserting dashes ``$\mn$'' in such a way that there is no position
in which both $\mu^1$ and $\mu^2$ have a dash. A \emph{match} is a
position where $\mu^1$ and $\mu^2$ have the same character, a
\emph{mismatch} is a position where $\mu^1$ and $\mu^2$ have
different characters, and a \emph{space} is a position in which one
of $\mu^1$ and $\mu^2$ has a dash. A simple scoring scheme consists
of two parameters $\alpha$ and $\beta$ denoting mismatch and space
penalties respectively. The reward of a match is set to $1$. The
score of an alignment with $z$ matches, $x$ mismatches, and $y$
spaces is then $z-x\alpha-y\beta$. Observe that these numbers always
satisfy $2z+2x+y=n_1+n_2$.

This model for sequence alignment can be translated into a probabilistic model, and is a particular case of a
so-called pair hidden Markov model. The problem of determining the
highest scoring alignment for given values of $\alpha$ and $\beta$
is equivalent to the inference problem in the pair hidden Markov
model, with some parameters set to functions of $\alpha$ and
$\beta$, or to $0$ or $1$. In this setting, an observation is a pair
of sequences $\mathbf{\tau}=(\sigma^1,\sigma^2)$, and the number of
observed variables is $n=n_1+n_2$. An explanation is then an optimal
alignment, since the values of the hidden variables indicate the
positions of the spaces.

In the rest of this chapter we will refer to this as the
\emph{$2$-parameter model for sequence alignment}. Note that it
actually comes from a $3$-parameter model where the reward for a
match has, without loss of generality, been set to $1$. The Newton
polytopes of the coordinates of the model are defined in a
$3$-dimensional space, but in fact they lie on a plane, as we will
see next. Thus, the parameter space has only two degrees of freedom.

For each pair of sequences $\mathbf{\tau}$, the Newton polytope of
the polynomial $f_\mathbf{\tau}$ is the convex hull of the points
$(x,y,z)$ whose coordinates are the number of mismatches, spaces,
and matches, respectively, of each possible alignment of the pair.
This polytope lies on the plane $2z+2x+y=n_1+n_2$, so no information
is lost by considering its projection onto the $xy$-plane instead.
This projection is just the convex hull of the points $(x,y)$ giving
the number of mismatches and spaces of each alignment. For any
alignment of sequences of lengths $n_1$ and $n_2$, the corresponding
point $(x,y)$ lies inside the square $[0,n]^2$, where $n=n_1+n_2$.
Therefore, since we are dealing with lattice polygons inside
$[0,n]^2$, it follows from Theorem~\ref{th:fif} that the number of
inference functions of this model is $O(n^2)$. Next we show that
this quadratic bound is tight, even in the case of the binary
alphabet.

\begin{prop}\label{prop:tight_binary}
Consider the $2$-parameter model for sequence alignment for two
observed sequences of length $n$ and let $\Sigma'=\{0,1\}$ be the
binary alphabet. Then, the number of inference functions of this
model is $\Theta(n^2)$.
\end{prop}

\begin{proof}
The above argument shows that $cn^2$ is an upper bound on the
number of inference functions of the model, for some constant $c$. To prove the
proposition, we will argue that there is some constant $c'$ such that there are at least $c'n^2$
such functions.

Since the two sequences have the same length, the number of spaces
in any alignment is even. For convenience, we define $y'=y/2$ and
$\beta'=2\beta$, and we will work with the coordinates $(x,y',z)$
and the parameters $\alpha$ and $\beta'$. The value $y'$ is called
the number of insertions (half the number of spaces), and $\beta'$
is the insertion penalty. For fixed values of $\alpha$ and $\beta'$,
the explanation of an observation
$\mathbf{\tau}=(\sigma^1,\sigma^2)$ is given by the vertex of
$\np(f_\mathbf{\tau})$ that is maximal in the direction of the
vector $(-\alpha,-\beta',1)$. In this model, $\np(f_\mathbf{\tau})$
is the convex hull of the points $(x,y',z)$ whose coordinates are
the number of mismatches, insertions and matches of the alignments
of $\sigma^1$ and $\sigma^2$.

The argument in the proof of Theorem~\ref{th:fif} shows that the
number of inference functions of this model is the number of cones
in the common refinement of the normal fans of
$\np(f_\mathbf{\tau})$, where $\mathbf{\tau}$ runs over all pairs of
sequences of length $n$ in the alphabet $\Sigma'$. Since the
polytopes $\np(f_\mathbf{\tau})$ lie on the plane $x+y'+z=n$, it is
equivalent to consider the normal fans of their projections onto the
$y'z$-plane. These projections are lattice polygons contained in the
square $[0,n]^2$. We denote by $P_\mathbf{\tau}$ the projection of
$\np(f_\mathbf{\tau})$ onto the $y'z$-plane.

We will construct a collection of pairs of binary sequences
$\mathbf{\tau}=(\sigma^1,\sigma^2)$ so that the total number of
different slopes of the edges of the polygons $\np(f_\mathbf{\tau})$
is $\Omega(n^2)$. This will imply that the number of cones in
$\bigwedge_{\mathbf{\tau}}\nf(\np(f_\mathbf{\tau}))$ is
$\Omega(n^2)$, where $\mathbf{\tau}$ ranges over all pairs of binary
sequences of length $n$.

We claim that for any positive integers $u$ and $v$ with $u<v$ and
$6v-2u\le n$, there exists a pair $\mathbf{\tau}$ of binary
sequences of length $n$ such that $P_\mathbf{\tau}$ has an edge of
slope $u/v$. This will imply that the number of different slopes
created by the edges of the polygons $P_\mathbf{\tau}$ is
$\Omega(n^2)$.

Thus, it only remains to prove the claim. Given positive integers
$u$ and $v$ as above, let $a:=2v$, $b:=v-u$. Assume first that
$n=6v-2u=2a+2b$. Consider the sequences
$$\sigma^1=0^a1^b0^b1^a,\qquad \sigma^2=1^a0^b1^b0^a,$$ where $0^a$
indicates that the symbol $0$ is repeated $a$ times. Let
$\mathbf{\tau}=(\sigma^1,\sigma^2)$. Then, it is not hard to see
that the polygon $P_\mathbf{\tau}$ for this pair of sequences has
four vertices: $v_0=(0,0)$, $v_1=(b,3b)$, $v_2=(a+b,a+b)$ and
$v_3=(n,0)$. The slope of the edge between $v_1$ and $v_2$ is
$(a-2b)/a=u/v$.

If $n>6v-2u=2a+2b$, we just append $0^{n-2a-2b}$ to both sequences
$\sigma^1$ and $\sigma^2$. In this case, the vertices of
$P_\mathbf{\tau}$ are $(0,n-2a-2b)$, $(b,n-2a+b)$, $(a+b,n-a-b)$,
$(n,0)$ and $(n-2a-2b,0)$.

\begin{figure}[h]
\begin{center} \includegraphics[width=0.35\textwidth]{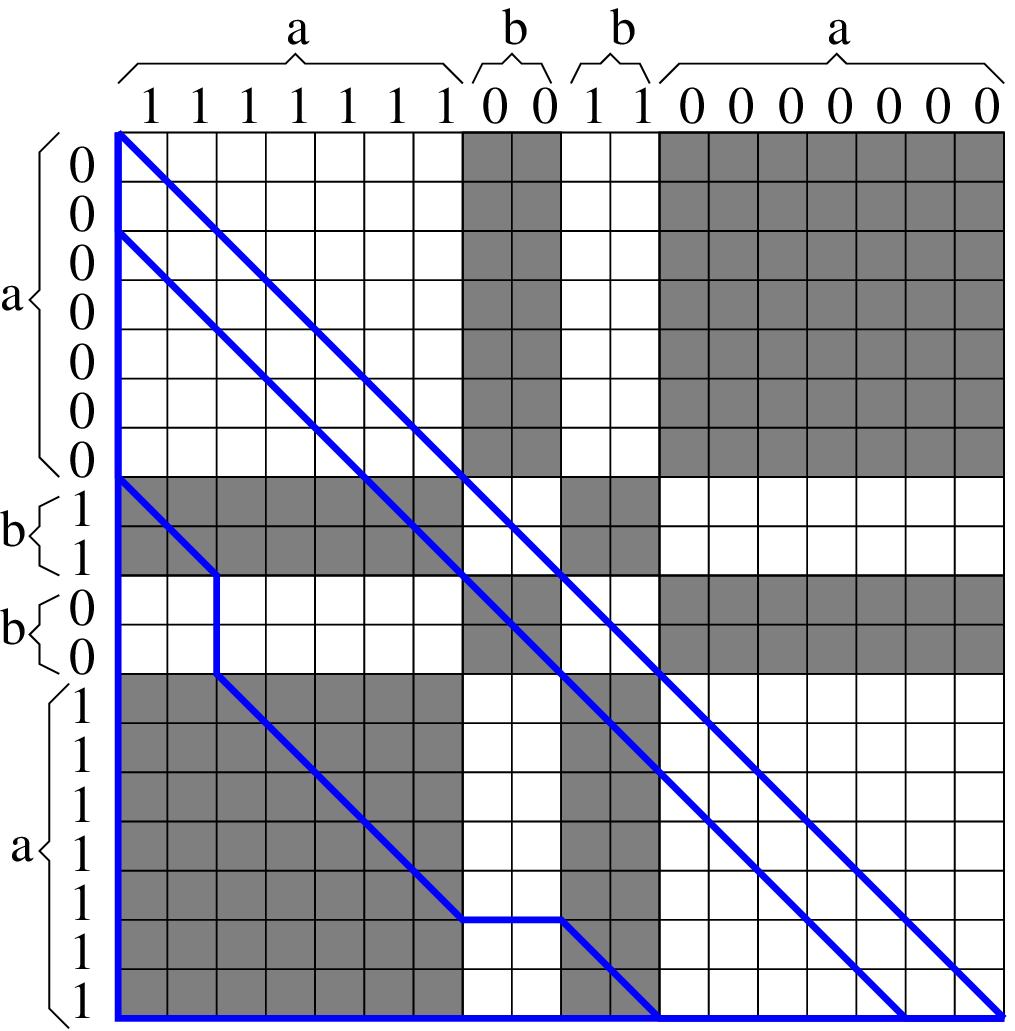}\qquad
\includegraphics[width=0.35\textwidth]{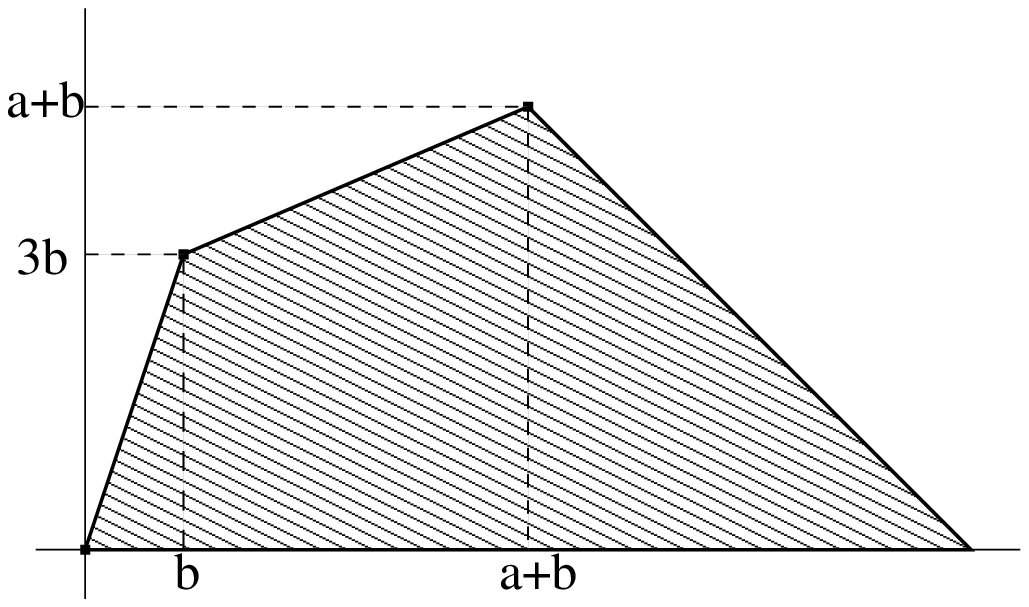}
\caption{A pair of binary sequences of length 18 giving the slope
$3/7$ in their alignment polytope. \label{fig:binary_slope}}
\end{center}
\end{figure}

Note that if $v-u$ is even, the construction can be done with
sequences of length $n=3v-u$ by taking $a:=v$, $b:=(v-u)/2$.
Figure~\ref{fig:binary_slope} shows the alignment graph and the
polygon $P_\mathbf{\tau}$ for $a=7$, $b=2$.
\end{proof}

In most cases, one is interested only in those inference functions
that are biologically meaningful. In our case, meaningful values of
the parameters occur when $\alpha,\beta\ge0$, which means that
mismatches and spaces are penalized instead of rewarded. Sometimes
one also requires that $\alpha\le\beta$, which means that a mismatch
should be penalized less than two spaces. It is interesting to
observe that our construction in the proof of
Proposition~\ref{prop:tight_binary} not only shows that the total
number of inference functions is $\Omega(n^2)$, but also that the
number of biologically meaningful ones is still $\Omega(n^2)$. This
is because the different rays created in our construction have a
biologically meaningful direction in the parameter space.

\section{Final remarks} \label{sec:final_remarks}

An interpretation of Theorem~\ref{th:fif} is that the ability to
change the values of the parameters of a graphical model does not
give as much freedom as it may appear. There is a very large number
of possible ways to assign an explanation to each observation.
However, only a tiny proportion of these come from a consistent
method for choosing the most probable explanation for a certain
choice of parameters. Even though the parameters can vary
continuously, the number of different inference functions that can
be obtained is at most polynomial in the number of edges of the
model, assuming that the number of parameters is fixed.

In the case of sequence alignment, the number of possible functions
that associate an alignment to each pair of sequences of length $n$
is doubly-exponential in $n$. However, the number of functions that
pick the alignment with highest score in the $2$-parameter model,
for some choice of the parameters $\alpha$ and $\beta$, is only
$\Theta(n^2)$. Thus, most ways of assigning alignments to pairs of
sequences do not correspond to any consistent choice of parameters.
If we use a model with more parameters, say $d$, the number of
inference functions may be larger, but still polynomial in $n$,
namely $O(n^{d(d-1)})$.

Having shown that the number of inference functions of a graphical
model is polynomial in the size of the model, an interesting next
step would be to find an efficient way to precompute all the
inference functions for given models. This would allow us to give
the answer (the explanation) to a query (an observation) very
quickly. It follows from this chapter that it is computationally
feasible to precompute the polytope $\np(\f)$, whose vertices
correspond to the inference functions. However, the difficulty
arises when we try to describe a particular inference function
efficiently. The problem is that the characterization of an
inference function involves an exponential number of observations.


\subsection*{Acknowledgements}

The authors are grateful to Graham Denham, Lior Pachter, Carl
Pomerance, Bernd Sturmfels, and Ravi Kannan for helpful discussions.
The first author was partially supported by the J. William Fulbright
Association of Spanish Fulbright Alumni.

\end{document}